\documentclass[a4paper, 11pt] {article}
\usepackage[cp1251]{inputenc}
\usepackage[russian, english]{babel}
\usepackage{amsfonts}
\usepackage{mathtext}
\usepackage{cite}

\usepackage{graphicx}

\usepackage{verbatim}
\usepackage{amsmath}
\usepackage{amsthm}
\usepackage{amssymb}
\usepackage{delarray}
\usepackage{mathrsfs}
\usepackage{xcolor}

\begin{document}


\begin{center}
{\bf\large ON NON-UNIQUENESS OF RECOVERING STURM--LIOUVILLE OPERATORS WITH DELAY AND THE NEUMANN BOUNDARY
CONDITION AT ZERO}
\end{center}

\begin{center}
{\bf\large Neboj\v{s}a Djuri\'c\footnote{Faculty of Architecture, Civil Engineering and Geodesy, University
of Banja Luka, {\it nebojsa.djuric@aggf.unibl.org}} and Sergey Buterin\footnote{Department of Mathematics,
Saratov State University, {\it buterinsa@info.sgu.ru}}}
\end{center}

{\bf Abstract.} As is known, for each fixed $\nu\in\{0,1\},$ the spectra of two operators generated by
$-y''(x)+q(x)y(x-a)$ and the boundary conditions $y^{(\nu)}(0)=y^{(j)}(\pi)=0,$ $j=0,1,$ uniquely determine
the complex-valued square-integrable potential $q(x)$ vanishing on $(0,a)$ as soon as $a\in[2\pi/5,\pi).$
Meanwhile, it actually became the {\it main question} of the inverse spectral theory for Sturm--Liouville
operators with constant delay whether the uniqueness holds also for smaller values of $a.$ Recently, a
negative answer was given by the authors [Appl. Math. Lett. 113 (2021) 106862] for $a\in[\pi/3,2\pi/5)$ in
the case $\nu=0$ by constructing an infinite family of iso-bispectral potentials. Moreover, an essential and
dramatic reason was established why this strategy, generally speaking, fails in the remarkable case when
$\nu=1.$ Here we construct a counterexample giving a {\it negative} answer for $\nu=1,$ which is an important
subcase of the Robin boundary condition at zero. We also refine the former counterexample for $\nu=0$ to
$W_2^1$-potentials.


\\

{\large\bf 1. Introduction}
\\

In recent years, there appeared an interest in inverse spectral problems for Sturm--Liouville-type operators
with deviating argument, see, e.g., papers \cite{Pik91, FrYur12, VladPik16, ButYur19, ButPikYur17, Ign18,
BondYur18-1, BMSh, VPV19, Yur19, DV19, Yur20, SS19, Dur, Yang, Wang, BK, Wang21, DB21}, among which a big
part is devoted to operators with delay. For $j=0,1,$ denote by $\{\lambda_{n,j}\}_{n\ge1}$ the spectrum of
the boundary value problem
\begin{equation}\label{1}
-y''(x)+q(x)y(x-a)=\lambda y(x),\quad 0<x<\pi,
\end{equation}
\begin{equation}\label{2}
U(y)=y^{(j)}(\pi)=0
\end{equation}
with delay $a\in(0,\pi)$ and a complex-valued potential $q(x)\in L_2(0,\pi)$ such that $q(x)=0$ on $(0,a),$
while $U(y)=y(0)$ or $U(y)=y'(0)-hy(0),$ $h\in{\mathbb C}.$ Such cases of $U(y)$ correspond to Dirichlet and
Robin boundary conditions at zero, to which we will refer as {\it Case 1} and case {\it Case 2},
respectively.

\smallskip
{\bf Inverse Problem 1.} Given the spectra $\{\lambda_{n,0}\}_{n\ge1}$ and $\{\lambda_{n,1}\}_{n\ge1},$ find
the potential $q(x).$

\smallskip
Alternatively, one can consider the case of Robin boundary conditions also at the right end:
\begin{equation}\label{rob}
U(y)=y'(\pi)+H_jy(\pi)=0,\quad j=0,1, \quad H_0,H_1\in{\mathbb C}, \quad  H_0\ne H_1,
\end{equation}
which, however, can be easily reduced to conditions (\ref{2}). Moreover, in the cases of Robin boundary
conditions, the coefficients $h$ and $H_0,$ $H_1$ are uniquely determined by the two spectra (see
\cite{Yur20}).

Various aspects of Inverse Problem~1 were studied in \cite{Pik91, FrYur12, VladPik16, ButYur19, ButPikYur17,
BondYur18-1, BMSh, VPV19, Yur19, DV19, Yur20, Dur, DB21} and other works. In particular, it is well known
that the two spectra uniquely determine the potential as soon as $a\ge\pi/2.$ Moreover, the inverse problem
is overdetermined (see \cite{ButYur19}). For $a<\pi/2,$ the dependence of the characteristic function of any
problem of the form (\ref{1}) and (\ref{2}) on the potential is nonlinear. It became actually the {\it main
question} of the inverse spectral theory for the operators with constant delay whether the uniqueness holds
also for small~$a.$ Recently, a positive answer when $a\in[2\pi/5,\pi/2)$ was given in \cite{BondYur18-1} for
Case~1 and independently in \cite{VPV19} for Case~2. However, recent authors' paper \cite{DB21} gave a
negative answer in Case~1 as soon as $a\in[\pi/3,2\pi/5).$ Specifically, for each such $a,$ we constructed an
infinite family of different iso-bispectral potentials~$q(x),$ i.e. for which both problems consisting of
(\ref{1}) and (\ref{2}) possess one and the same pair of spectra. This appeared quite unexpected because of
the inconsistence with Borg's classical uniqueness result for $a=0$ \cite{B}, and also in light of recent
paper \cite{Yur20} announcing the uniqueness for $a\in[\pi/3,\pi)$ in Case~2 for boundary conditions
(\ref{rob}). In \cite{DB21}, we also established an essential and dramatic reason why the idea of that
counterexample, generally speaking, fails in Case~2 (see Remark~2 in \cite{DB21}).

In the present paper, we return to Case~2 and construct a counterexample giving a {\it negative} answer for
the Neumann boundary condition at zero (i.e. when $h=0),$ which, unfortunately, refutes the uniqueness
theorem in \cite{Yur20} for $a\in[\pi/3,2\pi/5).$ For this purpose, we establish Theorem~1 (see the next
section) first, which reduces finding a counterexample in Case~2 to constructing a Hermitian integral
operator of a special form possessing an eigenfunction with the mean value zero. Even though the existence of
such an operator was highly believable, finding its concrete example appeared to be a quite difficult task.
After a series of computational experiments we constructed several numerical examples, one of which
fortunately admitted a precise elementary implementation (see Proposition~1).

This new non-uniqueness result along with the one in \cite{DB21} changes the further strategy of studying
inverse problems for the operators with delay (see also Remark~1 in \cite{DB21}). In particular, there
appears the relevance of finding various conditions on the class of potentials that would guarantee the
uniqueness of recovering $q(x).$ That is especially important for justifying constructive procedures for
solving Inverse Problem~1 when $a<2\pi/5,$ otherwise the corresponding algorithms in \cite{Yur19} and
\cite{Yur20} become indefinite. By virtue of (\ref{3.5}), as such a condition one can impose holomorphy of
$q(x)$ on an appropriate part of $(a,\pi).$

Finally, we note that both our counterexamples involve discontinuous potentials. So it is also relevant to
investigate the possibility of constructing iso-bispectral potentials in $W_2^k[0,\pi]$ with $k\in{\mathbb
N}$ so large as possible. Here we construct such potentials in Case~1 for $k=1$ and $a\in(\pi/3,2\pi/5)$ (see
Theorem~2).
\\

{\large\bf 2. The main results}
\\

For $\nu,j\in\{0,1\},$ denote by ${\cal L}_{\nu,j}(a,q)$ the eigenvalue problem for equation (\ref{1}) under
the boundary conditions
\begin{equation}\label{DirNeu}
y^{(\nu)}(0)=y^{(j)}(\pi)=0.
\end{equation}

Fix $a\in[\pi/3,2\pi/5).$ Following the main idea of the work \cite{DB21}, we consider the integral operator
\begin{equation}\label{2.1}
M_hf(x)=\int_\frac{3a}2^{\pi-x+\frac{a}2} K_h\Big(x+t-\frac{a}2\Big)f(t)\,dt, \;\; \frac{3a}2<x<\pi-a, \quad
{\rm where} \quad K_h(x)=\int_x^\pi h(\tau)\,d\tau,
\end{equation}
with a nonzero real-valued function $h(x)\in L_2(5a/2,\pi).$ Thus, $M_h$ is a nonzero compact Hermitian
operator in $L_2(3a/2,\pi-a)$ and, hence, it has at least one nonzero eigenvalue $\eta.$

Further, fix $\nu\in\{0,1\}$ and put $h_\nu(x):=(-1)^\nu h(x)/\eta.$ Then $(-1)^\nu$ is an eigenvalue of the
operator $M_{h_\nu}.$ Let $e_\nu(x)$ be some related eigenfunction,~i.e.
\begin{equation}\label{2.2}
M_{h_\nu}e_\nu(x)=(-1)^\nu e_\nu(x), \quad \frac{3a}2<x<\pi-a.
\end{equation}
Consider the one-parametric family of potentials $B_\nu:=\{q_{\alpha,\nu}(x)\}_{\alpha\in{\mathbb C}}$
determined by the formula
\begin{equation}\label{2.3}
q_{\alpha,\nu}(x)=\left\{\begin{array}{cl}\displaystyle 0,
 &\displaystyle x\in\Big(0,\frac{3a}2\Big)\cup(\pi-a,2a)\cup\Big(\pi-\frac{a}2,\frac{5a}2\Big),\\[3mm]
\displaystyle \alpha e_\nu(x), &\displaystyle x\in\Big(\frac{3a}2,\pi-a\Big),\\[3mm]
\displaystyle -\alpha K_{h_\nu}\Big(x+\frac{a}2\Big)\int_\frac{3a}2^{x-\frac{a}2} e_\nu(t)\,dt,
 &\displaystyle x\in\Big(2a,\pi-\frac{a}2\Big),\\[3mm]
h_\nu(x), &\displaystyle x\in\Big(\frac{5a}2,\pi\Big).
\end{array}\right.
\end{equation}
In \cite{DB21}, it was established that, for $j=0,1,$ the spectrum of the problem ${\cal
L}_{0,j}(a,q_{\alpha,0})$ is independent of~$\alpha$ for any function $h(x)$ conditioned above. This means
that $B_0$ is an iso-bispectral set (of potentials) for these two problems, i.e. the solution of Inverse
Problem~1 in Case~1 is not unique.

Moreover, in \cite{DB21}, it was noted that acting in the analogous way but for the problems ${\cal
L}_{1,j}(a,q),$ $j=0,1,$ would lead to the family of potentials $B_1.$ However, an essential reason was
established why $B_1,$ generally speaking, does not form an iso-bispectral set for these two problems (see
also Remark~1 in Section~3). In this paper, we find a concrete example when it is. We begin with the
following theorem.

\smallskip
{\bf Theorem 1. }{\it For $j=0,1,$ the spectrum of the problem ${\cal L}_{1,j}(a,q_{\alpha,1})$ is
independent of $\alpha$ as soon~as
\begin{equation}\label{2.3.1}
\int_\frac{3a}2^{\pi-a} e_1(x)\,dx=0.
\end{equation}
}

Thus, the problem of constructing iso-bispectral potentials for the problems ${\cal L}_{1,0}(a,q)$ and ${\cal
L}_{1,1}(a,q)$ is reduced to the question of finding a function $h(x)\in L_2(5a/2,\pi)$ such that the
operator $M_h$ has at least one eigenfunction possessing the zero mean value but related to a nonzero
eigenvalue. The answer to this question is given by the following assertion.

\smallskip
{\bf Proposition 1. }{\it Put
\begin{equation}\label{2.3.2}
h_1(x):=\frac{6\pi^2}{(2\pi-5a)^2} \cos\frac{\pi\sqrt{10}(\pi-x)}{2\pi-5a}, \quad
e_1(x):=\cos\frac{2\pi(2x-3a)}{2\pi-5a}+\cos\frac{\pi(2x-3a)}{2\pi-5a}.
\end{equation}
Then relation (\ref{2.2}) for $\nu=1$ as well as equality (\ref{2.3.1}) are fulfilled.}

\smallskip
Theorem~1 and Proposition~1 imply that the family $B_1$ constructed by using the functions $h_1(x)$ and
$e_1(x)$ that are determined by (\ref{2.3.2}) consists of iso-bispectral potentials for the problems ${\cal
L}_{1,0}(a,q)$ and ${\cal L}_{1,1}(a,q).$ Thus, Inverse Problem~1 is {\it not uniquely solvable also in
Case~2}.

Finally, we construct iso-bispectral potentials in $W_2^1[0,\pi]$ in Case~1. For this purpose, we consider
\begin{equation}\label{2.3.2.1}
e_0(x):=\sin\frac{2\pi(2x-3a)}{2\pi-5a}+2\sin\frac{\pi(2x-3a)}{2\pi-5a}
\end{equation}
and introduce the family of potentials $\tilde B_0:=\{q_\alpha(x)\}_{\alpha\in{\mathbb C}}$ determined by the
formula
\begin{equation}\label{2.3.3}
q_\alpha(x)=\left\{\begin{array}{cl}\displaystyle 0,
 &\displaystyle x\in\Big[0,\frac{3a}2\Big)\cup[\pi-a,2a),\\[3mm]
\displaystyle \alpha e_0(x), &\displaystyle x\in\Big[\frac{3a}2,\pi-a\Big),\\[3mm]
\displaystyle -\alpha K_{h_0}\Big(x+\frac{a}2\Big)\int_\frac{3a}2^{x-\frac{a}2} e_0(t)\,dt,
 &\displaystyle x\in\Big[2a,\pi-\frac{a}2\Big),\\[3mm]
\displaystyle g(x), &\displaystyle x\in\Big[\pi-\frac{a}2,\frac{5a}2\Big),\\[3mm]
h_0(x), &\displaystyle x\in\Big[\frac{5a}2,\pi\Big],
\end{array}\right.
\end{equation}
with $h_0(x)=h_1(x),$ where $h_1(x)$ is, in turn, determined in (\ref{2.3.2}), while $g(x)$ is an arbitrary
fixed function in $W_2^1[\pi-a/2,5a/2]$ obeying the boundary conditions
\begin{equation}\label{2.3.4}
g\Big(\pi-\frac{a}2\Big)=0, \quad g\Big(\frac{5a}2\Big)=\frac{6\pi^2}{(2\pi-5a)^2}\cos\frac{\pi\sqrt{10}}2.
\end{equation}

{\bf Theorem 2. }{\it For $j=0,1,$ the spectrum of the boundary value problem ${\cal L}_{0,j}(a,q_\alpha)$ is
independent of~$\alpha.$ Moreover, the inclusion $\tilde B_0\subset W_2^1[0,\pi]$ holds as soon as
$a\in(\pi/3,2\pi/5).$}

\smallskip
The proofs of Theorems~1 and~2 as well as Propositions~1 are provided in the next section.
\\

{\large\bf 3. The proofs}
\\

Denote by $y_0(x,\lambda)$ and $y_1(x,\lambda)$ the unique solutions of equation (\ref{1}) under the initial
conditions $y_\nu^{(j)}(0,\lambda)=\delta_{\nu,j},$ $\nu,j=0,1,$ where $\delta_{\nu,j}$ is the Kronecker
delta. For any pair of $\nu,j\in\{0,1\},$ eigenvalues of the boundary value problem ${\cal L}_{\nu,j}(a,q)$
with account of multiplicity coincide with zeros of the entire function
$\Delta_{\nu,j}(\lambda)=y_{1-\nu}^{(j)}(\pi,\lambda),$ which is called {\it characteristic function} of the
problem ${\cal L}_{\nu,j}(a,q).$ Thus, the spectrum of any problem ${\cal L}_{\nu,j}(a,q)$ does not depend on
$q(x)\in B$ for some subset $B\subset L_2(0,\pi)$ as soon as neither does the corresponding characteristic
function $\Delta_{\nu,j}(\lambda).$ Put $\rho^2=\lambda$ and denote
\begin{equation}\label{3.0}
\omega:=\int_a^\pi q(x)\,dx.
\end{equation}

Before proving Theorem~1, we provide necessary information from \cite{DB21}. For $\nu,j=0,1,$ we have
\begin{equation}\label{3.1}
\left.\begin{array}{c}
\displaystyle \Delta_{\nu,\nu}(\lambda)=(-\lambda)^{\nu}\Big(\frac{\sin\rho\pi}\rho-
\omega\frac{\cos\rho(\pi-a)}{2\lambda}
+\frac{(-1)^\nu}{2\lambda}\int_a^\pi w_\nu(x)\cos\rho(\pi-2x+a)\,dx\Big),\\[4mm]
\displaystyle \Delta_{\nu,j}(\lambda)=\cos\rho\pi +\omega\frac{\sin\rho(\pi-a)}{2\rho}
+\frac{(-1)^j}{2\rho}\int_a^\pi w_\nu(x)\sin\rho(\pi-2x+a)\,dx, \quad \nu\ne j, 
\end{array}\right\}
\end{equation}
where the functions $w_\nu(x)$ are determined by the formula
\begin{equation}\label{3.5}
w_\nu(x)=\left\{\begin{array}{l}
\displaystyle q(x),\quad x\in\Big(a,\frac{3a}2\Big)\cup\Big(\pi-\frac{a}2,\pi\Big),\\[3mm]
\displaystyle q(x)+Q_\nu(x),\quad x\in\Big(\frac{3a}2,\pi-\frac{a}2\Big),
\end{array}\right.
\end{equation}
while
\begin{equation}\label{3.6}
Q_\nu(x)=\int_a^{x-\frac{a}2} q(t)\,dt\int_{x+\frac{a}2}^\pi q(\tau)\,d\tau -(-1)^\nu\int_a^{\pi-x+\frac{a}2}
q(t)\,dt\int_{x+t-\frac{a}2}^\pi q(\tau)\,d\tau.
\end{equation}

\smallskip
{\bf Remark 1.} As was established \cite{DB21}, the difference between the cases $\nu=0$ and $\nu=1$ is as
folows. Since the functions $\Delta_{\nu,j}(\lambda)$ are entire in $\lambda,$ the first representation in
(\ref{3.1}) for $\nu=0$ implies
\begin{equation}\label{3.6.1}
\omega=\int_a^\pi w_0(x)\,dx,
\end{equation}
which can also be checked directly using (\ref{3.5}) and (\ref{3.6}) for $\nu=0.$ Thus, for $\nu=0,$ the
iso-bispectrality of~$B$ requires only~$w_0(x)$'s independence of $q(x)\in B.$ However, for $\nu=1,$ there is
no relation analogous to~(\ref{3.6.1}). In other words, the constant $\omega$ is not determined by $w_1(x).$
Thus, both functions $\Delta_{1,0}(\lambda)$ and $\Delta_{1,1}(\lambda)$ may depend on $q(x)\in B$ even when
$w_1(x)$ does not.

\smallskip
Let $q(x)=0$ on $(a,3a/2).$ Hence, formulae (\ref{3.5}) and (\ref{3.6}) give
\begin{equation}\label{3.9}
w_\nu(x)=\left\{\begin{array}{cl} \displaystyle 0,
& \displaystyle x\in\Big(a,\frac{3a}2\Big),\\[3mm]
 \displaystyle q(x)-(-1)^\nu M_hq(x),
& \displaystyle x\in\Big(\frac{3a}2,\pi-a\Big),\\[3mm]
\displaystyle q(x), & x\in(\pi-a,2a),\\[3mm]
\displaystyle q(x)+K_h\Big(x+\frac{a}2\Big)\int_\frac{3a}2^{x-\frac{a}2} q(t)\,dt,
 & \displaystyle x\in\Big(2a,\pi-\frac{a}2\Big),\\[3mm]
\displaystyle q(x),
 & \displaystyle x\in\Big(\pi-\frac{a}2,\frac{5a}2\Big),\\[3mm]
 \displaystyle h(x),  & \displaystyle x\in\Big(\frac{5a}2,\pi\Big),
\end{array}\right.
\end{equation}
where $h(x)=q(x)$ on $(5a/2,\pi),$ while $M_h$ and $K_h(x)$ are determined by (\ref{2.1}).

\smallskip
{\it Proof of Theorem~1.} Substituting $q(x)=q_{\alpha,1}(x)$ into~(\ref{3.9}) for $\nu=1,$ where
$q_{\alpha,1}(x)$ is determined by (\ref{2.3}) with $\nu=1,$ and taking (\ref{2.2}) for $\nu=1$ into account,
we arrive at
$$
w_1(x)=0, \quad a<x<\frac{5a}2, \qquad w_1(x)=h_1(x),\quad \frac{5a}2<x<\pi.
$$
Thus, the function $w_1(x)$ is independent of $\alpha.$ Hence, it remains to prove that so is also the
value~$\omega$ determined by formula (\ref{3.0}) with $q(x)=q_{\alpha,1}(x).$ Integrating the third line in
(\ref{2.3}) for $\nu=1,$ we get
$$
{\cal I}:=\int_{2a}^{\pi-\frac{a}2}K_{h_1}\Big(x+\frac{a}2\Big)\,dx\int_\frac{3a}2^{x-\frac{a}2} e_1(t)\,dt
=\int_{2a}^{\pi-\frac{a}2}K_{h_1}\Big(x+\frac{a}2\Big)\,dx \int_\frac{3a}2^{x-\frac{a}2} e_1(x+a-t)\,dt.
$$
Changing the order of integration and then the internal integration variable, we calculate
$$
{\cal I}=\int_\frac{3a}2^{\pi-\frac{a}2}dx \int_{x+\frac{a}2}^{\pi-\frac{a}2}K_{h_1}\Big(t+\frac{a}2\Big)
e_1(t+a-x)\,dt = \int_\frac{3a}2^{\pi-a}dx \int_\frac{3a}2^{\pi-x+\frac{a}2}K_{h_1}\Big(x+t-\frac{a}2\Big)
e_1(t)\,dt,
$$
which along with the first equality in (\ref{2.1}) as well as (\ref{2.2}) for $\nu=1$ and (\ref{2.3.1})
implies
$$
{\cal I}= \int_\frac{3a}2^{\pi-a} M_{h_1} e_1(x)\,dx =-\int_\frac{3a}2^{\pi-a} e_1(x)\,dx=0.
$$
Thus, according to (\ref{2.3}) for $\nu=1,$ the assumption (\ref{2.3.1}) of the theorem gives
$$
\omega=\int_a^\pi q_{\alpha,1}(x)\,dx=\int_\frac{5a}2^\pi h_1(x)\,dx,
$$
i.e. the value $\omega$ does not depend on $\alpha,$ which finishes the proof. $\hfill\Box$

\smallskip
For shortening the remaining proofs, we provide the following auxiliary assertion.

\smallskip
{\bf Proposition 2. }{\it For each fixed $\nu\in\{0,1\},$ relation (\ref{2.2}) is equivalent to the relation
\begin{equation}\label{2.2.1}
m_{\chi_\nu}\epsilon_\nu(\xi)=(-1)^\nu \epsilon_\nu(\xi), \quad 0<\xi<1, \quad m_\chi f(\xi):=\int_0^{1-\xi}
f(\eta)\,d\eta\int_{\xi+\eta}^1 \chi(\theta)\,d\theta,
\end{equation}
as soon as
\begin{equation}\label{2.1.1}
\epsilon_\nu(\xi)=e_\nu\Big(\frac{3a}2 +\Big(\pi-\frac{5a}2\Big)\xi\Big), \quad
\chi_\nu(\theta)=\Big(\pi-\frac{5a}2\Big)^2h_\nu\Big(\frac{5a}2 +\Big(\pi-\frac{5a}2\Big)\theta\Big).
\end{equation}
}

\smallskip
{\it Proof.} Making in (\ref{2.2}) the change of variable $\xi:=(2x-3a)/A,$ where $A=2\pi-5a,$ we obtain
$$
M_{h_\nu}e_\nu\Big(\frac{3a}2 +\Big(\pi-\frac{5a}2\Big)\xi\Big)=(-1)^\nu\epsilon_\nu(\xi), \quad 0<\xi<1.
$$
Using (\ref{2.1}) and successively making the changes $\eta:=(2t-3a)/A$ and $\theta:=(2\tau-5a)/A,$ we get
$$
(-1)^\nu\epsilon_\nu(\xi)= \Big(\pi-\frac{5a}2\Big)\int_0^{1-\xi}\epsilon_\nu(\eta)\,d\eta
\int_{\frac{5a}2+(\pi-\frac{5a}2)(\xi+\eta)}^\pi h_\nu(\tau)\,d\tau =\int_0^{1-\xi}\epsilon_\nu(\eta)\,d\eta
\int_{\xi+\eta}^1 \chi_\nu(\theta)\,d\theta,
$$
which coincides with (\ref{2.2.1}). $\hfill\Box$

\smallskip
{\it Proof of Proposition~1.} Let us start with (\ref{2.3.1}), which can be checked by the direct
substitution:
$$
\int_\frac{3a}2^{\pi-a} e_1(x)\,dx=\frac{2\pi-5a}{2\pi}\Big(\frac12
\sin\frac{2\pi(2x-3a)}{2\pi-5a}+\sin\frac{\pi(2x-3a)}{2\pi-5a}\Big)\Big|_{x=\frac{3a}2}^{\pi-a}=0.
$$
According to Proposition~2, it remains to prove relation (\ref{2.2.1}) for $\nu=1$ with the functions
$$
\epsilon_1(\xi)=\cos\pi\xi+\cos2\pi\xi, \quad
\chi_1(\theta)=\frac{3\pi^2}2\cos\frac{\pi\sqrt{10}(1-\theta)}2,
$$
which are determined by (\ref{2.3.2}) and (\ref{2.1.1}) for $\nu=1.$ Indeed, it is easy to calculate
$$
m_{\chi_1}\epsilon(\xi)=\frac{3\pi}{2\sqrt{10}}(A_1+A_2), \quad A_j=2\int_0^{1-\xi}\cos\pi j\eta
\cdot\sin\frac{\pi\sqrt{10}(1-\xi-\eta)}2\, d\eta
$$
$$
=\sum_{k=0}^1\int_0^{1-\xi}\sin\Big(\frac{\pi\sqrt{10}(1-\xi)}2 -\pi\Big(\frac{\sqrt{10}}2+(-1)^k
j\Big)\eta\Big)\, d\eta = \frac1\pi \frac{2\sqrt{10}}{5-2j^2} \Big(\cos\pi j(1-\xi)
-\cos\frac{\pi\sqrt{10}(1-\xi)}2\Big),
$$
where $j=1,2,$ which leads to (\ref{2.2.1}) for $\nu=1.$ $\hfill\Box$

\smallskip
{\it Proof of Theorem~2.} Using Proposition~2 as in the preceding proof, one can establish~(\ref{2.2}) for
$\nu=0$ with $h_0(x)=h_1(x)$ determined in (\ref{2.3.2}) and $e_0(x)$ determined by (\ref{2.3.2.1}). Further,
substituting $q(x)=q_\alpha(x)$ into~(\ref{3.9}) for $\nu=0,$ where $q_\alpha(x)$ is determined by
(\ref{2.3.3}), and taking~(\ref{2.2}) for $\nu=0$ into account, we get
$$
w_0(x)=0, \;\; a<x<\pi-\frac{a}2, \quad w_0(x)= g(x), \;\; \pi-\frac{a}2<x<\frac{5a}2, \quad w_0(x)=h_0(x),
\;\; \frac{5a}2<x<\pi.
$$
Thus, the function $w_0(x)$ is independent of $\alpha.$ Moreover, according to (\ref{3.1}) and (\ref{3.6.1}),
for $j=0,1,$ the characteristic function $\Delta_{0,j}(\lambda)$ of the problem ${\cal L}_{0,j}(a,q_\alpha)$
is independent of $\alpha.$

Finally, we note that, for any $\alpha\in{\mathbb C}$ and $a\in(\pi/3,2\pi/5),$ the inclusion $q_\alpha(x)\in
W_2^1[0,\pi]$ follows from (\ref{2.3.2})--(\ref{2.3.4}) along with the last equality in (\ref{2.1}).
$\hfill\Box$
\\

{\bf Acknowledgement.} The first author was supported by the Project 19.032/961-103/19 of the Republic of
Srpska Ministry for Scientific and Technological Development, Higher Education and Information Society. The
second author was supported by Grants 19-01-00102 and 20-31-70005 of the Russian Foundation for Basic
Research.

\end{document}